\documentclass[preprints,article,accept,oneauthor,pdftex,10pt,a4paper]{mdpi} 

\firstpage{1} 
\pubvolume{xx}
\issuenum{1}
\articlenumber{5}
\pubyear{2018}
\copyrightyear{2018}
\history{Received: date; Accepted: date; Published: date}

\Title{On the relation between topological entropy and restoration entropy}

\Author{Christoph Kawan $^{1}$}

\AuthorNames{Christoph Kawan}

\address{%
$^{1}$ \quad Universit\"{a}t Passau, Fakult\"{a}t f\"{u}r Informatik und Mathematik, Innstra{\ss}e 33, 94032 Passau (Germany); christoph.kawan@uni-passau.de}

\corres{Correspondence: christoph.kawan@uni-passau.de; Tel.: +49(0)851 509-3363}

\abstract{In the context of state estimation under communication constraints, several notions of dynamical entropy play a fundamental role, among them: topological entropy and restoration entropy. In this paper, we present a theorem which demonstrates that for most dynamical systems restoration entropy strictly exceeds topological entropy. This implies that robust estimation policies in general require a higher rate of data transmission than non-robust ones. The proof of our theorem is quite short, but uses sophisticated tools from the theory of smooth dynamical systems.}


\newcommand{\tp}{\mathrm{top}}%
\newcommand{\rst}{\mathrm{res}}%
\newcommand{\id}{\mathrm{id}}%
\newcommand{\cl}{\mathrm{cl}}%
\newcommand{\inner}{\mathrm{int}}%
\newcommand{\md}{\mathrm{mod}}%

\newcommand{\rmd}{\mathrm{d}}%
\newcommand{\rmD}{\mathrm{D}}%
\newcommand{\rme}{\mathrm{e}}%

\newcommand{\N}{\mathbb{N}}%
\newcommand{\Z}{\mathbb{Z}}%
\newcommand{\R}{\mathbb{R}}%
\newcommand{\T}{\mathbb{T}}%

\newcommand{\ep}{\varepsilon}%
\newcommand{\tm}{\times}%

\newcommand{\CC}{\mathcal{C}}%
\newcommand{\EC}{\mathcal{E}}%
\newcommand{\FC}{\mathcal{F}}%
\newcommand{\MC}{\mathcal{M}}%
\newcommand{\PC}{\mathcal{P}}%
\newcommand{\QC}{\mathcal{Q}}%

\begin{document}

\section{Introduction}

This paper compares two notions of entropy that are relevant in the context of state estimation under communication constraints. Since the work of Savkin \cite{Sav}, it is well-known that the topological entropy of a dynamical system characterizes the smallest rate of information above which an estimator, receiving its state information at this rate, is able to generate a state estimate of arbitrary precision. Topological entropy is a quantity that has been studied in the mathematical field of dynamical systems since the 1960s and has turned out to be a useful tool for solving many theoretical and practical problems, cf.~the survey \cite{Kat} and the monograph \cite{Dow}. A big drawback of this notion in the context of state estimation is that topological entropy is highly discontinuous with respect to the dynamical system under consideration in any reasonable topology, cf.~\cite{Mis}. As a consequence, estimation policies based on topological entropy are likely to suffer from a lack of robustness. Additionally, topological entropy is very hard to compute or estimate. There are only few numerical approaches that potentially work for multi-dimensional systems, cf.~\cite{COH,Dea,FJO,NPi}, and each of them has its drawbacks and restrictions.%

A possible remedy for these problems is provided in the works \cite{MP1,MP2} of Matveev and Pogromsky. One of the main ideas in these papers is to replace the topological entropy as a figure-of-merit for the necessary rate of data transmission with a possibly larger quantity, named restoration entropy, which describes the smallest data rate above which a more robust form of state estimation can be achieved (called \emph{regular observability} in \cite{MP1,MP2}).%

Looking at one of the simplest types of nonlinear dynamical systems, namely Anosov diffeomorphisms, the main result of the paper at hand demonstrates that for most dynamical systems we have to expect that the restoration entropy strictly exceeds the topological entropy. That is, to achieve a state estimation objective that is more robust with respect to perturbations, one has to pay the price of using a channel that allows for a larger rate of data transmission. More specifically, our result shows that the equality of topological and restoration entropy implies a great amount of uniformity in the dynamical system under consideration, which can be expressed in terms of the unstable Lyapunov exponents at each point, whose sum essentially has to be a constant. Such a property can easily be destroyed by a small perturbation, showing that arbitrarily close to the given system we find systems whose restoration entropy strictly exceed their topological entropy.%

To prove our result, we need a number of high-level concepts and results from the theory of topological, measurable and smooth dynamical systems. This includes the concepts of topological and metric pressure, Lyapunov exponents, SRB measures and uniform hyperbolicity.%

For further reading on the topic of state estimation under communication constraints, we refer the reader to \cite{KYu,LMi,Mat,MP1,MP2,MSa,Sav} and the references given therein.%

The structure of this paper is as follows: In Section \ref{sec_prelim}, we collect all necessary definitions and results from the theory of dynamical systems. Section \ref{sec_resent} introduces the concept of restoration entropy and explains its operational meaning in the context of estimation under communication constraints. In Section \ref{sec_results}, we prove our main result and provide some interpretation and an example. Finally, Section \ref{sec_discuss} contains some concluding remarks.%

\section{Tools from dynamical systems}\label{sec_prelim}

{\bf Notation:} By $\Z$ we denote the set of all integers, by $\N$ the set of positive integers and $\N_0 := \{0\} \cup \N$. All logarithms are taken to the base $2$. If $M$ is a Riemannian manifold, we write $|\cdot|$ for the induced norm on any tangent space $T_xM$, $x\in M$. The notation $\|\cdot\|$ is reserved for operator norms. We write $\cl A$ and $\inner A$ for the closure and the interior of a set $A$ in a metric space, respectively.%
\vspace{0.4cm}

In this paper, we use several sophisticated results from the theory of dynamical systems, in particular from smooth ergodic theory. In the following, we try to explain these results without going too much into technical details.%

Let $T:X \rightarrow X$ be a continuous map on a compact metric space $(X,d)$. Via its iterates%
\begin{equation*}
  T^0 := \id_X,\quad T^{n+1} := T \circ T^n,\quad n = 0,1,2,\ldots%
\end{equation*}
the map $T$ generates a discrete-time dynamical system on $X$ with associated orbits $\{T^n(x)\}_{n\in\N_0}$, $x\in X$. We call the pair $(X,T)$ a \emph{topological dynamical system}, briefly a \emph{TDS}.%

\subsection{Entropy and pressure}

Let $(X,T)$ be a TDS. The \emph{topological entropy} $h_{\tp}(T)$ measures the total exponential complexity of the orbit structure of $(X,T)$ in terms of the maximal numbers of finite-time orbits that are distinguishable w.r.t.~to a finite resolution. One amongst different possible formal definitions is as follows. For $n\in\N$ and $\ep>0$, a set $E \subset X$ is called $(n,\ep,T)$-separated if for any $x,y\in E$ with $x \neq y$ we have%
\begin{equation*}
  d(T^i(x),T^i(y)) \geq \ep \mbox{\quad for at least one\ } 0 \leq i < n.%
\end{equation*}
That is, we can distinguish any two points in $E$ at a resolution of $\ep$ by looking at their length-$n$ finite-time orbits. By compactness of $X$, there is a uniform upper bound on the cardinality of any $(n,\ep,T)$-separated set. Writing $r(n,\ep,T)$ for the maximal possible cardinality,%
\begin{equation*}
  h_{\tp}(T) := \lim_{\ep\downarrow0}\limsup_{n \rightarrow \infty}\frac{1}{n}\log r(n,\ep,T).%
\end{equation*}
This definition is due to Bowen \cite{Bow} and (independently) Dinaburg \cite{Din}. However, it should be noted that the first definition of topological entropy, given by Adler, Konheim and McAndrew \cite{AKM}, was in terms of open covers of $X$ and was modeled in strict analogy to the metric ($=$ measure-theoretic) entropy defined earlier by Kolmogorov and Sinai \cite{Kol,Sin}.%

To define metric entropy, one additionally needs a Borel probability measure $\mu$ on $X$ that is preserved by $T$ in the sense that $\mu(A) = \mu(T^{-1}(A))$ for every Borel set $A$. By the theorem of Krylov-Bogolyubov, every continuous map on a compact space admits at least one such measure, cf.~\cite[Thm.~4.1.1]{KHa}. We write $\MC_T$ for the set of all $T$-invariant Borel probability measures. For any finite measurable partition $\PC$ of $X$, we define the entropy of $T$ on $\PC$ by%
\begin{equation*}
  h_{\mu}(T;\PC) := \lim_{n \rightarrow \infty}\frac{1}{n}H_{\mu}\Bigl(\bigvee_{i=0}^{n-1}T^{-i}\PC\Bigr).%
\end{equation*}
Here $\bigvee$ denotes the join operation. That is, $\bigvee_{i=0}^{n-1}T^{-i}\PC$ is the partition of $X$ whose elements are all intersections of the form $P_0 \cap T^{-1}(P_1) \cap \ldots \cap T^{-n+1}(P_{n-1})$ with $P_i \in \PC$. Moreover, $H_{\mu}(\cdot)$ denotes the Shannon entropy of a partition, i.e., $H_{\mu}(\QC) = -\sum_{Q\in\QC}\mu(Q)\log\mu(Q)$ for any finite partition $\QC$. The \emph{metric entropy} of $T$ w.r.t.~$\mu$ is then defined by%
\begin{equation*}
  h_{\mu}(T) := \sup_{\PC}h_{\mu}(T;\PC),%
\end{equation*}
the supremum taken over all finite measurable partitions $\PC$ of $X$.\footnote{Replacing measurable partitions with open covers and Shannon entropy with the logarithm of the cardinality of a minimal finite subcover, the same construction yields the topological entropy as defined in \cite{AKM}.}%

To understand the meaning of $h_{\mu}$, note that $H_{\mu}(\QC)$ is the average amount of uncertainty as one attempts to predict the partition element to which a randomly chosen point belongs. Hence, $h_{\mu}(T)$ measures the average uncertainty per iteration in guessing the partition element of a typical length $n$-orbit.%

The \emph{variational principle for entropy} states that%
\begin{equation}\label{eq_entropy_vp}
  h_{\tp}(T) = \sup_{\mu\in\MC_T}h_{\mu}(T),%
\end{equation}
where the supremum is not necessarily a maximum. This variational principle can be regarded as a quantitative version of the theorem of Krylov-Bogolyubov.%

Another concept (of which entropy is a special case) used in dynamical systems and inspired by ideas in thermodynamics is \emph{pressure}. In this context, any continuous function $\phi:X \rightarrow \R$, also called a \emph{potential} or an \emph{observable}, gives rises to the \emph{metric pressure} of $T$ w.r.t.~$\phi$ for a given $\mu \in \MC_T$, defined as%
\begin{equation*}
  P_{\mu}(T,\phi) := h_{\mu}(T) + \int\phi\rmd\mu.%
\end{equation*}
To define an associated notion of \emph{topological pressure}, put $S_n\phi(x) := \sum_{i=0}^{n-1}\phi(T^i(x))$ and%
\begin{equation*}
  R(n,\ep,\phi;T) := \sup\Bigl\{ \sum_{x\in E}\rme^{S_n\phi(x)} : E \subset X \mbox{ is } (n,\ep,T)\mbox{-separated} \Bigr\}.%
\end{equation*}
Then the topological pressure of $T$ w.r.t.~$\phi$ is given by%
\begin{equation*}
  P_{\tp}(T,\phi) := \lim_{\ep\downarrow0}\limsup_{n \rightarrow \infty}\frac{1}{n}\log R(n,\ep,\phi;T).%
\end{equation*}
The associated variational principle, first proved in \cite{Wal}, reads%
\begin{equation}\label{eq_pressure_vp}
  P_{\tp}(T,\phi) = \sup_{\mu \in \MC_T}P_{\mu}(T,\phi),%
\end{equation}
which includes \eqref{eq_entropy_vp} a special case (simply put $\phi = 0$).%

\subsection{Subadditive cocycles}

Let $T:X \rightarrow X$ be a map. A \emph{subadditive cocycle} over $(X,T)$ is a sequence $(f_n)_{n\in\N_0}$ of functions $f_n:X \rightarrow \R$ satisfying%
\begin{equation*}
  f_{n+m}(x) \leq f_n(x) + f_m(T^n(x)),\quad \forall n,m\in\N_0,\ x \in X.%
\end{equation*}
If equality holds in this relation, we call $(f_n)_{n\in\N_0}$ an \emph{additive cocycle} over $(X,T)$.%

If $X$ has the structure of a probability space with a $\sigma$-algebra $\FC$ and a probability measure $\mu$ on $\FC$, $T$ is measurable and $\mu$ is $T$-invariant, we speak of a \emph{measurable subadditive cocycle} provided that all $f_n$ are measurable. In the context of a TDS $(X,T)$, we speak of a \emph{continuous subadditive cocycle} if all $f_n$ are continuous.%

The most fundamental result about subadditive cocycles is Kingman's Subadditive Ergodic Theorem, cf.~\cite[Thm.~2.1.4]{Dow}:%

\begin{Theorem}
Let $T:X \rightarrow X$ be a measure-preserving map on a probability space $(X,\FC,\mu)$ and $(f_n)_{n\in\N_0}$ a measurable subadditive cocycle over $(X,T)$ such that each $f_n$ is integrable. Then the limit%
\begin{equation*}
  \lim_{n \rightarrow \infty}\frac{1}{n}f_n(x)%
\end{equation*}
exists for $\mu$-almost every $x\in X$. If, additionally, $\mu$ is ergodic, then the limit is constant with%
\begin{equation}\label{eq_kset}
  \lim_{n \rightarrow \infty}\frac{1}{n}f_n(x) = \lim_{n \rightarrow \infty}\frac{1}{n}\int f_n \rmd\mu.%
\end{equation}
\end{Theorem}

Observe that the limit on the right-hand side of \eqref{eq_kset} always exists by Fekete's subadditivity lemma (see \cite[Fact 2.1.1]{Dow}), because the sequence $a_n := \int f_n \rmd\mu$ is subadditive, i.e., $a_{n+m} \leq a_n + a_m$. Kingman's Theorem can, in particular, be applied if $(X,T)$ is a TDS, $\mu \in \MC_T$ and $(f_n)_{n\in\N_0}$ is a continuous subadditive cocycle.%

Now we consider again a TDS $(X,T)$ and a continuous subadditive cocycle $(f_n)_{n\in\N_0}$ over $(X,T)$. We define the extremal growth rate of $(f_n)$ by%
\begin{equation*}
  \beta[(f_n)] := \sup_{x\in X}\limsup_{n \rightarrow \infty}\frac{1}{n}f_n(x).%
\end{equation*}
The following result is well-known, and can be found in \cite[Thm.~A.3]{Mor}, for instance:%

\begin{Lemma}\label{lem_sc_ergodic}
Let $(f_n)_{n\in\N_0}$ be a continuous subadditive coycle over a TDS $(X,T)$. Then%
\begin{equation*}
  \beta[(f_n)] = \sup_{\mu\in\MC_T}\inf_{n>0}\frac{1}{n}\int f_n \rmd\mu = \inf_{n>0}\sup_{x\in X}\frac{1}{n}f_n(x) = \inf_{n>0}\sup_{\mu \in \MC_T}\frac{1}{n}\int f_n \rmd\mu.%
\end{equation*}
Here, all infima can be replaced with limits. Moreover, every supremum is attained.%
\end{Lemma}

\subsection{Lyapunov exponents, SRB measures and Pesin's formula}

To describe the long-term dynamical behavior of smooth systems, the notion of \emph{Lyapunov exponents} is crucial. Given a $C^1$-diffeomorphism $T:M \rightarrow M$ on a compact Riemannian manifold $M$, the \emph{Lyapunov exponent} at $x \in M$ in direction $0 \neq v \in T_xM$ is the number%
\begin{equation*}
  \lambda(x,v) := \lim_{n \rightarrow \infty}\frac{1}{n}\log|\rmD T^n(x)v|,%
\end{equation*}
provided that the limit exists. Lyapunov exponents measure how fast nearby solutions diverge from each other. The most general result on their existence and their properties is the \emph{Multiplicative Ergodic Theorem (MET)}, also known as \emph{Oseledets Theorem}, cf.~\cite{Arn,CKl}. We need the following version of the theorem (which is not the most general):%

\begin{Theorem}\label{thm_MET}
Let $T:M \rightarrow M$ be a $C^1$-diffeomorphism of a compact Riemannian manifold $M$ and $\mu \in \MC_T$. Then there exists a Borel set $\Omega \subset M$ with $\mu(\Omega)=1$ and $T(\Omega) = \Omega$ such that the following holds: For every $x\in\Omega$, there exist numbers $\lambda_1(x) > \ldots > \lambda_{r(x)}(x)$ and the tangent space at $x$ splits into linear subspaces as%
\begin{equation*}
  T_xM = E_1(x) \oplus \cdots \oplus E_{r(x)}(x)%
\end{equation*}
such that the following properties hold:%
\begin{enumerate}
\item[(i)] For every $0 \neq v \in E_i(x)$ we have%
\begin{equation*}
  \lim_{n \rightarrow \pm \infty}\frac{1}{n}\log|\rmD T^n(x)v| = \lambda_i(x).%
\end{equation*}
\item[(ii)] The functions $r(\cdot)$, $\dim E_i(\cdot)$ and $\lambda_i(\cdot)$ are measurable and constant along orbits. Moreover,%
\begin{equation*}
  \rmD T(x)E_i(x) = E_i(T(x)),\quad i=1,\ldots,r(x).%
\end{equation*}
\item[(iii)] For every $x\in\Omega$ the limit%
\begin{equation*}
  \Lambda_x := \lim_{n \rightarrow \infty}(\rmD T^n(x)^* \rmD T^n(x))^{1/2n}%
\end{equation*}
exists and the different eigenvalues of $\Lambda_x$ are $2^{\lambda_1(x)},\ldots,2^{\lambda_{r(x)}(x)}$.%
\end{enumerate}
\end{Theorem}

Typically, a given map has a huge number of associated invariant measures. To obtain a good description of the global dynamical behavior, one has to select specific invariant measures that determine the behavior of the system on a large set of initial states. In this context, the notion of an \emph{SRB measure (Sinai-Ruelle-Bowen measure)} comes into play. An SRB measure is a measure with at least one positive Lyapunov exponent almost everywhere, having absolutely continuous conditional measures on unstable manifolds. We are not going to give a technical definition of the latter property. Instead, we state the following celebrated theorem due to Ledrappier and Young \cite{LYo}, which characterizes this property in terms of metric entropy. Here we use the short-cut%
\begin{equation*}
  \lambda^+(x) := \sum_{i=1}^{r(x)} \max\{0,\lambda_i(x)\dim E_i(x)\}%
\end{equation*}
for the sum of all positive Lyapunov exponents at a point $x\in\Omega$, counted with multiplicities.%

\begin{Theorem}\label{thm_ledrappier_young}
Let $T:M \rightarrow M$ be a $C^2$-diffeomorphism of a compact manifold $M$ and $\mu \in \MC_T$. Then the formula%
\begin{equation}\label{eq_pesin}
  h_{\mu}(T) = \int \lambda^+ \rmd\mu%
\end{equation}
holds if and only if $\mu$ has absolutely continuous conditional measures on unstable manifolds.%
\end{Theorem}

Additionally, note that for any $C^1$-diffeomorphism $T$ and any $\mu\in\MC_T$, the inequality%
\begin{equation}\label{eq_ruelle}
  h_{\mu}(T) \leq \int \lambda^+ \rmd\mu%
\end{equation}
holds, which is known as \emph{Ruelle's inequality} or \emph{Ruelle-Margulis inequality} \cite{Rue}. (The formula \eqref{eq_pesin} was first proved by Pesin for smooth invariant measures.)%

\subsection{Anosov diffeomorphisms}

One of the simplest classes of smooth dynamical systems with complicated dynamical behavior is the class of Anosov diffeomorphisms. In this paper, we use these systems for two reasons. First, they have positive topological entropy, and second they are very well understood and a lot of tools are available to describe their properties.%

Let $M$ be a compact Riemannian manifold. A $C^1$-diffeomorphism $T:M \rightarrow M$ is called an \emph{Anosov diffeomorphism} if there exists a splitting%
\begin{equation*}
  T_xM = E^u_x \oplus E^s_x,\quad \forall x \in M%
\end{equation*}
into linear subspaces such that the following conditions are satisfied:%
\begin{enumerate}
\item[(1)] $\rmD T(x)E^u_x = E^u_{T(x)}$ and $\rmD T(x)E^s_x = E^s_{T(x)}$ for all $x\in M$.%
\item[(2)] There are constants $c \geq 1$ and $\lambda \in (0,1)$ so that for all $x\in M$ and $n\in\N_0$,%
\begin{align*}
  |\rmD T^n(x)v| &\leq c \lambda^n |v| \mbox{ for all } v \in E^s_x,\\
	|\rmD T^{-n}(x)v| &\leq c \lambda^n |v| \mbox{ for all } v \in E^u_x.%
\end{align*}
\end{enumerate}
From (1) and (2) it automatically follows that $E^s_x$ and $E^u_x$ vary continuously with $x$, cf.~\cite[Prop.~6.4.4]{KHa}. The existence of a splitting as above is also known as \emph{uniform hyperbolicity}.%

The simplest examples of Anosov diffeomorphisms are hyperbolic linear torus automorphisms, i.e., maps on the $n$-dimensional torus $\T^n = \R^n/\Z^n$ of the form%
\begin{equation*}
  T_A(x) = Ax\ \ (\md\ \Z^n),\quad T_A:\T^n \rightarrow \T^n,%
\end{equation*}
where $A \in \Z^{n\tm n}$ is an integer matrix satisfying $|\det A| = 1$ and $|\lambda| \neq 1$ for all eigenvalues $\lambda$ of $A$. Observe that the assumption $|\det A|=1$ guarantees that $T_A$ is invertible with inverse $T_A^{-1} = T_{A^{-1}}$ (because $A^{-1}$ also has integer entries) and at the same time implies that $T_A$ is area-preserving. That is, the normalized Lebesgue measure on $\T^n$ is an element of $\MC_{T_A}$. The assumption on the eigenvalues of $A$ together with the fact that the derivative $\rmD T_A(x)$ at any point $x\in \T^n$ can be identified with $A$ itself implies the Anosov properties (1) and (2).%

It is well-known that Anosov diffeomorphisms are structurally stable, i.e., any sufficiently small $C^1$-perturbation $T_{\ep}$ of an Anosov diffeomorphism $T:M \rightarrow M$ is also an Anosov diffeomorphism which is \emph{topologically conjugate} to $T$, see \cite[Prop.~6.4.6 and Cor.~18.2.2]{KHa}. That is, there exists a homeomorphism $h:M \rightarrow M$ so that%
\begin{equation*}
  h^{-1} \circ T_{\ep} \circ h = T.%
\end{equation*}

If we assume that $T$ is an arbitrary Anosov diffeomorphism of the torus, the existence of a unique entropy-maximizing measure $\mu$ follows. That is, $\mu$ is the unique element of $\MC_T$ satisfying%
\begin{equation*}
  h_{\tp}(T) = h_{\mu}(T).%
\end{equation*}
This follows from a combination of results that can be found in Katok \& Hasselblatt \cite{KHa}, namely Theorem 20.3.7, Proposition 18.6.5, Theorem 18.3.9 and Corollary 6.4.10. The entropy-maximizing measure $\mu$ is also known as the \emph{Bowen-measure}.%

In this context, also the notion of \emph{topological mixing} is important. An Anosov diffeomorphism (or simply a continuous map) $T:M \rightarrow M$ is called topologically mixing if for any two nonempty open sets $A,B \subset M$ there exists an integer $N$ such that $T^n(A) \cap B \neq \emptyset$ for all $n \geq N$. In particular, all Anosov diffeomorphisms on $\T^n$ are topologically mixing \cite[Prop.~18.6.5]{KHa}.%

\section{State estimation and restoration entropy}\label{sec_resent}

The notion of restoration entropy is introduced in \cite{MP2} for systems given by ODEs on $\R^n$. However, it is immediately clear from the definition that restoration entropy can be defined for any continuous map on a compact metric space as follows. Let $T:X \rightarrow X$ be a continuous map on a metric space $(X,d)$ and $K \subset X$ a compact set with $T(K)\subset K$. For every $x\in X$, $n\in\N$ and $\ep>0$, let $p(n,x,\ep)$ denote the smallest number of $\ep$-balls needed to cover the image $T^n(B_{\ep}(x) \cap K)$. If the map is not clear from the context, we also write $p(n,x,\ep;T)$. Then%
\begin{equation*}
  h_{\rst}(T_{|K}) := \lim_{n\rightarrow\infty}\frac{1}{n}\limsup_{\ep\downarrow0}\sup_{x\in X}\log p(n,x,\ep).%
\end{equation*}
The existence of the limit in $n$ follows from a subadditivity argument. If we assume that $T$ is a $C^1$-diffeomorphism of a compact Riemannian manifold, the numbers $p(n,x,\ep)$ can be estimated in terms of the unstable singular values of $\rmD T^n(x)$. This is related to the simple fact that the image of a ball under a linear map (in our case, the local linear approximation $\rmD T^n(x)$ to $T^n$) is an ellipsoid with semi-axes of lengths proportional to the singular values. This leads to the following result, proved in \cite[Thm.~11]{MP2} for continuous-time systems. The proofs carries over to discrete-time systems on Riemannian manifolds without any problem.%

\begin{Theorem}\label{thm_resent_form}
Let $T:M \rightarrow M$ be a $C^1$-diffeomorphism of a $d$-dimensional Riemannian manifold $M$ and $K \subset M$ a forward-invariant compact set of $T$ with $\cl K = \cl(\inner K)$. Then%
\begin{equation*}
  h_{\rst}(T_{|K}) = \lim_{n \rightarrow \infty}\frac{1}{n}\max_{x\in K}\sum_{i=1}^d \max\{0,\log\alpha_i(n,x)\},%
\end{equation*}
where $\alpha_1(n,x) \geq \ldots \geq \alpha_d(n,x)$ denote the singular values of $\rmD T^n(x)$.%
\end{Theorem}

For the analysis of $h_{\rst}$, based on the above formula, the following observations are crucial:%
\begin{itemize}
\item We have%
\begin{equation*}
  \sum_{i=1}^d \max\{0,\log \alpha_i(n,x)\} = \log \prod_{i=1}^d \max\{1,\log \alpha_i(n,x)\} = \log\|\rmD T^n(x)^{\wedge}\|,%
\end{equation*}
where $\rmD T^n(x)^{\wedge}$ denotes the linear map induced by $\rmD T^n(x)$ between the full exterior algebras of the tangent spaces $T_xM$ and $T_{T^n(x)}M$, respectively, see \cite[Ch.~I, Prop.~7.4.2]{BLR}.%
\item The sequence $f_n(x) := \log\|\rmD T^n(x)^{\wedge}\|$, $f_n:M \rightarrow \R$, is a continuous subadditive cocycle over $(K,T_{|K})$, since%
\begin{align*}
  f_{n+m}(x) &= \log\|\rmD T^{n+m}(x)^{\wedge}\| = \log\|\rmD T^m(T^n(x))^{\wedge} \rmD T^n(x)^{\wedge}\|\\
	           &\leq \log\left(\|\rmD T^m(T^n(x))^{\wedge}\| \cdot \|\rmD T^n(x)^{\wedge}\|\right)\\
	           &= \log\|\rmD T^n(x)^{\wedge}\| + \log\|\rmD T^m(T^n(x))^{\wedge}\| = f_n(x) + f_m(T^n(x)).%
\end{align*}
Alternatively, this follows from Horn's inequality for singular values, see \cite[Ch.~I, Prop.~2.3.1]{BLR}.%
\end{itemize}

In the following, we explain the operational meaning of the quantity $h_{\rst}(T_{|K})$. 

Consider the dynamical system given by%
\begin{equation}\label{eq_ds}
  x_{t+1} = T(x_t),\quad x_0\in K,\ t = 0,1,2,\ldots%
\end{equation}
Suppose that a sensor, fully observing the state $x_t$, sends its data to an encoder. At the sampling times $t=0,1,2,\ldots$, the encoder sends a signal $e_t$ through a noisefree discrete channel to a decoder (without transmission delay). The decoder acts as an observer of the system, trying to reconstruct the state from the received data. We write $\hat{x}_t$ for the estimate generated by the observer at time $t$. Moreover, we assume that we start with an initial estimate $\hat{x}_0 \in K$ of a specified accuracy.%

With $\MC$ denoting the coding alphabet, the encoder and the observer are described by mappings%
\begin{equation*}
  e_t = \CC_t(x_0,x_1,\ldots,x_t;\hat{x}_0,\delta),\quad \CC_t:K^{t+1} \tm K \tm \R_{>0} \rightarrow \MC,%
\end{equation*}
and%
\begin{equation*}
  \hat{x}_t = \EC_t(e_0,e_1,\ldots,e_t;\hat{x}_0,\delta),\quad \EC_t:\MC^{t+1} \tm K \tm \R_{>0} \rightarrow X.%
\end{equation*}
The argument $\delta$ corresponds to the initial error at time zero, i.e.\ $d(x_0,\hat{x}_0) \leq \delta$. In particular, we assume that both the encoder and the observer are given the data $\hat{x}_0$ and $\delta$.%

We assume that the channel can transmit at least $b_-(r)$ and at most $b_+(r)$ bits in any time interval of length $r$. The \emph{capacity} of the channel is then defined by%
\begin{equation*}
  C := \lim_{r\rightarrow\infty}\frac{b_-(r)}{r} = \lim_{r\rightarrow\infty}\frac{b_+(r)}{r},%
\end{equation*}
assuming that these limits exist and coincide.%

We consider the following two observation objectives:%
\begin{enumerate}
\item[(O1)] The observer \emph{observes} the system with exactness $\ep>0$ if there exists $\delta = \delta(\ep,K)$ so that $x_0,\hat{x}_0 \in K$ with $d(x_0,\hat{x}_0) \leq \delta$ implies%
\begin{equation*}
  \sup_{t\geq0}d(x_t,\hat{x}_t) \leq \ep.%
\end{equation*}
\item[(O2)] The observer \emph{regularly observes} the system if there exist $G,\delta_*>0$ so that for all $\delta \in (0,\delta_*)$ and $x_0,\hat{x}_0 \in K$ with $d(x_0,\hat{x}_0) \leq \delta$,%
\begin{equation*}
  \sup_{t\geq0}d(x_t,\hat{x}_t) \leq G\delta.%
\end{equation*}
\end{enumerate}

We say that the system is%
\begin{itemize}
\item \emph{observable on $K$} over a channel of capacity $C$ if for every $\ep>0$ an observer exists which observes the system with exactness $\ep$ over this channel;%
\item \emph{regularly observable on $K$} over a channel of capacity $C$ if there exists an observer which regularly observes the system over this channel.%
\end{itemize}

Then we have the following data-rate theorem, cf.~\cite[Thm.~8]{MP1} and \cite[Thm.~9]{MP2}.%

\begin{Theorem}
The smallest channel capacity $C_0$, so that system \eqref{eq_ds} is%
\begin{itemize}
\item observable on $K$ over every channel of capacity $C > C_0$ is given by%
\begin{equation*}
  C_0 = h_{\tp}(T_{|K}).%
\end{equation*}
\item regularly observable on $K$ over every channel of capacity $C > C_0$ is given by%
\begin{equation*}
  C_0 = h_{\rst}(T_{|K}).%
\end{equation*}
\end{itemize}
\end{Theorem}

Since regular observability implies observability, it is clear that%
\begin{equation*}
  h_{\tp}(T_{|K}) \leq h_{\rst}(T_{|K}).%
\end{equation*}

As already pointed out in the introduction, the quantity $h_{\tp}(\cdot)$ is highly discontinuous w.r.t.~the dynamical system. Moreover, the corresponding data-rate theorem has the disadvantage that the final error $\ep$ may be much larger than the initial error $\delta$, which cannot happen in the case of regular observability. From Theorem \ref{thm_resent_form} in combination with Lemma \ref{lem_sc_ergodic}, one sees that in the smooth case, $h_{\rst}$ is an infimum over functions that are continuous w.r.t.~$T$ in the $C^1$-topology. This implies at least upper semicontinuity. Hence, we can expect that coding and estimation strategies based on restoration entropy enjoy better properties than those based on topological entropy.%

\section{Results}\label{sec_results}

Before we present our main result, we prove two lemmas which are of independent interest.%

\begin{Lemma}\label{lem1}
Let $T:M \rightarrow M$ be a $C^2$-diffeomorphism on a compact Riemannian manifold $M$. Then for any $\mu\in\MC_T$ we have%
\begin{equation*}
  \int \lambda^+ \rmd\mu = \lim_{n \rightarrow \infty}\frac{1}{n}\int\log\|\rmD T^n(x)^{\wedge}\|\rmd\mu(x).%
\end{equation*}
\end{Lemma}

\begin{proof}
Let $d = \dim M$. First observe that we have the identity%
\begin{equation*}
  \|\rmD T^n(x)^{\wedge}\| = \max\Bigl\{1,\max_{1 \leq k \leq d} \prod_{i=1}^k \alpha_i(n,x)\Bigr\},%
\end{equation*}
where $\alpha_1(n,x) \geq \ldots \geq \alpha_d(n,x)$ are the singular values of $\rmD T^n(x)$, see \cite[Ch.~I, Prop.~7.4.2]{BLR}. Hence,%
\begin{equation*}
  \log\|\rmD T^n(x)^{\wedge}\| = \max\Bigl\{0,\max_{1 \leq k \leq d}\sum_{i=1}^k \log \alpha_i(n,x)\Bigr\}.%
\end{equation*}
The maximum over $k$ is clearly attained when $k$ is the maximal number such that $\alpha_i(n,x) > 1$ for all $1 \leq i \leq k$. Hence,%
\begin{equation*}
  \log\|\rmD T^n(x)^{\wedge}\| = \max\Bigl\{0,\sum_{\alpha_i(n,x)>1} \log \alpha_i(n,x)\Bigr\}.%
\end{equation*}
The numbers $\alpha_i(n,x)$ are the eigenvalues of $A_n(x) := (\rmD T^n(x)^* \rmD T^n(x))^{1/2}$. Theorem \ref{thm_MET} states that $A_n(x)^{1/n} \rightarrow \Lambda_x$ for $\mu$-almost every $x\in M$ and the logarithms of the eigenvalues of $\Lambda_x$ are the Lyapunov exponents at $x$. Since eigenvalues depend continuously on the matrix, it follows that%
\begin{equation*}
  \lim_{n \rightarrow \infty}\frac{1}{n}\log\|\rmD T^n(x)^{\wedge}\| = \lambda^+(x) \quad \mu\mbox{-a.e.}%
\end{equation*}
and consequently%
\begin{equation*}
  \int \lambda^+ \rmd\mu = \int\lim_{n \rightarrow \infty}\frac{1}{n}\log\|\rmD T^n(x)^{\wedge}\|\rmd\mu(x).%
\end{equation*}
Applying the Theorem of Dominated Convergence then yields the result.%
\end{proof}

\begin{Lemma}\label{lem2}
Let $T:M \rightarrow M$ be a $C^2$-diffeomorphism on a compact Riemannian manifold $M$ such that $h_{\tp}(T) = h_{\rst}(T)$. Then, if $T$ has an entropy-maximizing measure $\mu_*$, it follows that%
\begin{equation*}
  h_{\mu_*}(T) = \int \lambda^+ \rmd\mu_*.%
\end{equation*}
\end{Lemma}

\begin{proof}
Assume to the contrary that $h_{\mu_*}(T) < \int \lambda^+ \rmd\mu_*$ (using Ruelle's inequality \eqref{eq_ruelle}). Then Lemma \ref{lem1} implies%
\begin{equation*}
  h_{\tp}(T) = h_{\mu_*}(T) < \int \lambda^+ \rmd\mu_* = \lim_{n \rightarrow \infty}\frac{1}{n} \int \log\|\rmD T^n(x)^{\wedge}\| \rmd\mu_*(x).%
\end{equation*}
According to Theorem \ref{thm_resent_form} and the subsequent observation, an application of Lemma \ref{lem_sc_ergodic} yields%
\begin{equation*}
  h_{\rst}(T) = \sup_{\mu \in \MC_T}\lim_{n\rightarrow\infty}\frac{1}{n}\int \log \|\rmD T^n(x)^{\wedge}\| \rmd \mu(x).%
\end{equation*}
Combining these observations gives $h_{\tp}(T) < h_{\rst}(T)$, in contradiction to our assumption.%
\end{proof}

Now we are in position to state our main result.%

\begin{Theorem}\label{thm_mt}
Let $T:M \rightarrow M$ be a topologically mixing $C^2$-Anosov diffeomorphism on a compact Riemannian manifold $M$ such that $h_{\tp}(T) = h_{\rst}(T)$. Then the unique entropy-maximizing measure $\mu_* \in \MC_T$ is an SRB measure. Moreover, the function%
\begin{equation*}
  \mu \mapsto \int\lambda^+\rmd\mu,\quad \MC_T \rightarrow \R_{\geq0}%
\end{equation*}
is constant.%
\end{Theorem}

\begin{proof}
First note that the existence and uniqueness of an entropy-maximizing measure $\mu_*$ follows from \cite[Thm.~20.3.7, Thm.~18.3.9 and Cor.~6.4.10]{KHa}. Here the assumption that $T$ is topologically mixing is crucial. By the preceding lemma combined with Theorem \ref{thm_ledrappier_young} we already know that $\mu_*$ has absolutely continuous conditional measures on unstable manifolds. Since an Anosov diffeomorphism has positive Lyapunov exponents everywhere (where they exist), attained in all directions of the unstable subspace $E^u_x$, it follows that $\mu_*$ is an SRB measure.%

Now let $\mu \in \MC_T$ be chosen arbitrarily. Due to the invariance of $\mu$, we have%
\begin{align*}
  \int \log|\det \rmD T(x)_{|E^u_x}| \rmd \mu(x) &= \int \frac{1}{n} \sum_{i=0}^{n-1} \log |\det\rmD T(T^i(x))_{|E^u_{T^i(x)}}| \rmd\mu(x)\\
	&= \int \frac{1}{n}\log|\det\rmD T^n(x)_{|E^u_x}|\rmd \mu(x)%
\end{align*}
for every $n\in\N$, implying%
\begin{align*}
  \int \lambda^+ \rmd\mu &= \int \lim_{n \rightarrow \infty}\frac{1}{n}\log|\det\rmD T^n(x)_{|E^u_x}| \rmd\mu(x)\\
	&= \lim_{n\rightarrow\infty}\int \frac{1}{n}\log|\det\rmD T^n(x)_{|E^u_x}| \rmd\mu(x) = \int \log|\det \rmD T(x)_{|E^u_x}| \rmd \mu(x),%
\end{align*}
where we use Kingman's Subadditive Ergodic Theorem, applied to the continuous additive cocycle $f_n(x) := \log|\det\rmD T^n(x)_{|E^u_x}|$ ($n\in\N_0$), and the Theorem of Dominated Convergence. Observe that the function $J^uT(x) := \log|\det\rmD T(x)_{|E^u_x}|$ is continuous (using the fact that $x \mapsto E^u_x$ is continuous). Hence, we can consider the affine function%
\begin{equation*}
  \alpha_{\mu}:\R \rightarrow \R,\quad \alpha_{\mu}(t) := P_{\mu}(T,-t J^uT) = h_{\mu}(T) - t \int \lambda^+ \rmd \mu.%
\end{equation*}
The variational principle \eqref{eq_pressure_vp} for pressure tells us that%
\begin{equation}\label{eq_sup}
  P_{\tp}(-t J^uT) = \sup_{\mu \in \MC_T}\alpha_{\mu}(t),\quad \forall t\in\R.%
\end{equation}
Hence, $t \mapsto P_{\tp}(-t J^uT)$, as the supremum over affine functions, is a convex function.%

Using that $\mu_*$ is the entropy-maximizing measure and Theorem \ref{thm_ledrappier_young}, respectively, we obtain%
\begin{equation*}
  \alpha_{\mu_*}(0) = h_{\tp}(T) \mbox{\quad and \quad} \alpha_{\mu_*}(1) = 0.%
\end{equation*}
On the other hand, also%
\begin{equation*}
  P_{\tp}(- 0 \cdot \lambda^+) = h_{\tp}(T) \mbox{\quad and \quad} P_{\tp}(- 1 \cdot J^uT) = 0.%
\end{equation*}
The second identity here follows from the fact that $P_{\tp}(- 1 \cdot J^uT) = \sup_{\mu \in \MC_T}(h_{\mu}(T) - \int \lambda^+ \rmd\mu)$ and $h_{\mu}(T) \leq \int\lambda^+\rmd\mu$ by Ruelle's inequality \eqref{eq_ruelle}. Hence, $P_{\tp}(-1 \cdot J^uT) = h_{\mu_*}(T) - \int \lambda^+ \rmd\mu_* = 0$.%
 
By convexity of $t \mapsto P_{\tp}(-t J^uT)$ and \eqref{eq_sup}, this implies%
\begin{equation*}
  P_{\tp}(-tJ^uT) = \alpha_{\mu_*}(t),\quad \forall t\in\R.%
\end{equation*}
From \eqref{eq_sup} it now follows that all of the maps $\alpha_{\mu}$ have the same slope, i.e., $\int \lambda^+ \rmd\mu$ is independent of $\mu$.%
\end{proof}

The above theorem shows that the equality $h_{\tp}(T) = h_{\rst}(T)$ is a very restrictive condition. Indeed, this can be seen as follows. Any topologically mixing Anosov diffeomorphism has an abundance of periodic points. Indeed, the set of periodic points is dense in $M$, see \cite[Cor.~6.4.19]{KHa}. If we consider a periodic point $p \in M$ of period $n_p\in\N$, we can consider the invariant measure $\mu_p$ given by%
\begin{equation*}
  \mu_p := \frac{1}{n_p}\sum_{i=0}^{n_p-1} \delta_{T^i(p)}%
\end{equation*}
with $\delta_{(\cdot)}$ being the Dirac measure at a point. The above theorem implies that, under $h_{\tp}(T) = h_{\rst}(T)$, the number%
\begin{equation*}
  \gamma(p) := \int \lambda^+ \rmd\mu_p = \frac{1}{n_p} \log\left|\det \left( \rmD T^{n_p}(p)_{|E^u_p}:E^u_p \rightarrow E^u_p \right)\right|%
\end{equation*}
is independent of the periodic point $p$ chosen. On the other hand, we know that every sufficiently small $C^2$-perturbation of $T$ yields another $C^2$-Anosov diffeomorphism, topologically conjugate to $T$, hence also topologically mixing. If this perturbation is only performed in a small vicinity of a fixed periodic orbit, it can easily change the the number $\gamma(p)$, while not changing it for most of the other periodic orbits. As a consequence, the perturbed diffeomorphism $T_{\ep}$ cannot satisfy $h_{\tp}(T_{\ep}) = h_{\rst}(T_{\ep})$.%

The following corollary gives another characterization of Anosov diffeomorphisms with $h_{\tp} = h_{\rst}$ in a two-dimensional case.%

\begin{Corollary}\label{cor_c1conj}
Consider a $C^2$-area preserving Anosov diffeomorphism $T:\T^2 \rightarrow \T^2$ of the $2$-torus. Then the equality $h_{\tp}(T) = h_{\rst}(T)$ is equivalent to the existence of a hyperbolic linear automorphism $T_A:\T^2 \rightarrow \T^2$ and a $C^1$-diffeomorphism $h:\T^2 \rightarrow \T^2$ such that $h^{-1} \circ T \circ h = T_A$.%
\end{Corollary}

\begin{proof}
It follows immediately from Theorem \ref{thm_mt} in combination with \cite[Cor.~20.4.4]{KHa} that the identity $h_{\tp}(T) = h_{\rst}(T)$ implies the existence of a $C^1$-conjugacy as asserted. The other direction is easy to see, using the definition of restoration entropy. If $h^{-1} \circ T \circ h = T_A$, then also $h^{-1} \circ T^n \circ h = T_A^n$ for all $n\in\N$. We use that a $C^1$-map on a compact manifold has a global Lipschitz constant. Let $L := \mathrm{Lip}(h)$ and $L' := \mathrm{Lip}(h^{-1})$ be Lipschitz constants of $h$ and $h^{-1}$, respectively. Then%
\begin{equation*}
  T^n(B_{\ep}(x)) = h \circ T_A^n \circ h^{-1}(B_{\ep}(x)).%
\end{equation*}
Observe that $h^{-1}(B_{\ep}(x)) \subset B_{L'\ep}(h^{-1}(x))$. Let $N(l)$ denote the minimal number of $\ep$-balls needed to cover an $l\ep$-ball in $\T^2$ for any $l>0$. Then the minimal number of $\ep$-balls needed to cover $T_A^n h^{-1}(B_{\ep}(x))$ is bounded from above by $N(L') \max_{z\in \T^2}p(n,z,\ep;T_A)$. This implies%
\begin{equation*}
  p(n,x,\ep;T) \leq N(L)N(L') \max_{z\in\T^2}p(n,z,\ep;T_A).%
\end{equation*}
Hence,%
\begin{equation*}
  \sup_{x\in\T^2}\frac{1}{n}\log p(n,x,\ep;T) \leq \frac{1}{n}\log N(L)N(L') + \sup_{x\in\T^2}\frac{1}{n}\log p(n,x,\ep;T_A).%
\end{equation*}
Taking the $\limsup$ for $\ep\downarrow0$ and subsequently the limit for $n\rightarrow\infty$, we obtain that $h_{\rst}(T) \leq h_{\rst}(T_A)$. The other inequality can be proved analogously, so%
\begin{equation*}
  h_{\rst}(T) = h_{\rst}(T_A).%
\end{equation*}
Since $T$ and $T_A$ are topologically conjugate (the $C^1$-diffeomorphism $h$ is a homeomorphism, in particular), they also have the same topological entropy:%
\begin{equation*}
  h_{\tp}(T) = h_{\tp}(T_A).%
\end{equation*}
To complete the proof, it now suffices to show that $h_{\rst}(T_A) = h_{\tp}(T_A)$. We can compute $h_{\rst}(T_A)$ using Theorem \ref{thm_resent_form}. To this end, observe that $A$ is a hyperbolic matrix. If $|\lambda_1| > 1 > |\lambda_2|$ are its eigenvalues, we obtain%
\begin{equation*}
  \lim_{n \rightarrow \infty}\frac{1}{n}\sum_{i=1}^2 \max\{0,\log \alpha_i(n,x)\} = \log|\lambda_1| \quad \forall x \in \T^2,%
\end{equation*}
implying $h_{\rst}(T_A) = \log|\lambda_1|$. It is well-known that this is also the value of the topological entropy $h_{\tp}(T_A)$, see \cite[Sec.~4]{KHa}. This also follows from the combination of the variational principle with Theorem \ref{thm_ledrappier_young}.%
\end{proof}

The following example demonstrates how restrictive the condition $h_{\rst}(T) = h_{\tp}(T)$ is by looking at small perturbations of Arnold's Cat Map.%

\begin{Example}
Arnold's Cat Map is the hyperbolic linear $2$-torus automorphism $T_A:\T^2 \rightarrow \T^2$ induced by the integer matrix%
\begin{equation*}
  A := \left(\begin{array}{cc} 2 & 1 \\ 1 & 1 \end{array}\right)%
\end{equation*}
with determinant $\det A = 1$. Observe that the derivative $\rmD T_A(x)$ can be identified with $A$ for each $x \in \T^2$. Since $A$ is a hyperbolic matrix with eigenvalues%
\begin{equation*}
  \gamma_1 = \frac{3}{2} - \frac{1}{2}\sqrt{5} \mbox{\quad and\quad} \gamma_2 = \frac{3}{2} + \frac{1}{2}\sqrt{5}%
\end{equation*}
satisfying $|\gamma_2| > 1 > |\gamma_1|$, it follows that $T_A$ is a $C^{\infty}$-area preserving Anosov diffeomorphism. Hence, Corollary \ref{cor_c1conj} yields%
\begin{equation*}
  h_{\tp}(T_A) = h_{\rst}(T_A) = \log|\gamma_2|.%
\end{equation*}
Now we consider a perturbation of the form%
\begin{equation*}
  T_A^{\ep}(x,y) := \left(2x + y + \ep \sin(2\pi x),x + y\right)\ (\md\ \Z^2),\quad \ep > 0%
\end{equation*}
which is well-defined as a torus map, since the sine function is $2\pi$-periodic. By structural stability of Anosov diffeomorphisms, for a sufficiently small $\ep$, this map is topologically conjugate to $T_A$, hence has the same topological entropy $\log|\gamma_2|$. However, its restoration entropy is strictly greater. This can be seen by looking at the fixed point $(0,0)$ with associated derivative%
\begin{equation*}
  D T_A^{\ep}(0,0) = \left(\begin{array}{cc} 2 + 2\pi\ep & 1 \\ 1 & 1 \end{array}\right).%
\end{equation*}
The associated eigenvalues can be computed as%
\begin{equation*}
  \lambda_{\pm} = \frac{3}{2} + \pi\ep \pm \frac{1}{2}\sqrt{5 + 4\pi\ep(1 + \pi\ep)}.%
\end{equation*}
Since $\lambda_+ > \gamma_2$, Theorem \ref{thm_resent_form} yields $h_{\rst}(T_A^{\ep}) \geq \log|\lambda_+| > h_{\tp}(T_A^{\ep})$ for $\ep>0$ sufficiently small.%
\end{Example}

\section{Concluding remarks}\label{sec_discuss}

In this paper, we compared two notions of entropy for dynamical systems that have an operational meaning in the context of state estimation over digital channels: topological entropy and restoration entropy. Looking at Anosov diffeomorphisms, our main result demonstrates that the equality of these two quantities implies a great amount of uniformity in the given system. For area-preserving Anosov diffeomorphisms on the $2$-torus, this uniformity can be expressed in terms of the existence of a $C^1$-conjugacy to a linear system. Hence, we can conclude that for most dynamical systems the strict inequality $h_{\tp} < h_{\rst}$ holds. The operational meaning of this inequality is that for regular observability, as defined in Section \ref{sec_resent}, a strictly larger channel capacity is necessary than for observability.%

\vspace{6pt} 

\funding{This research received no external funding.}

\acknowledgments{The author owes particular thanks to Katrin Gelfert who provided one of the main ideas in the proof of Theorem \ref{thm_mt} during the Mini-Workshop \emph{Entropy, Information and Control} held at the Mathematisches Forschungsinstitut Oberwolfach from March 4 to March 10, 2018. The author also thanks Alexander Pogromsky for fruitful discussions on restoration entropy.}

\conflictsofinterest{The author declares no conflict of interest.}




\reftitle{References}




\end{document}